\documentstyle[12pt,amssymb,amscd]{amsart}  % amslatex

\newtheorem{thm}{Theorem}[section]
\newtheorem{cor}[thm]{Corollary}
\newtheorem{prop}[thm]{Proposition}
\newtheorem{lemma}[thm]{Lemma}

\theoremstyle{remark}
\newtheorem{remark}[thm]{Remark}
\newtheorem{example}[thm]{Example}

\theoremstyle{definition}
\newtheorem{defn}[thm]{Definition}

\numberwithin{equation}{section}

\newcommand{\bbC}{{\Bbb C}}

\newcommand{\bbL}{{\Bbb L}}
\newcommand{\bbR}{{\Bbb R}}

\newcommand{\bbZ}{{\Bbb Z}}

\newcommand{\bbP}{{\Bbb P}}

\newcommand{\cE}{{\cal E}}

\newcommand{\cO}{{\cal O}}

\newcommand{\cL}{{\cal L}}
\newcommand{\cM}{{\cal M}}

\newcommand{\cA}{{\cal A}}
\newcommand{\cB}{{\cal B}}

\newcommand{\cC}{{\cal C}}

\newcommand{\cK}{{\cal K}}

\newcommand{\supp}{\operatorname{Supp}}

\newcommand{\Hom}{\operatorname{Hom}}
\newcommand{\Ext}{\operatorname{Ext}}

\newcommand{\Perf}{\operatorname{Perf}}

\newcommand{\printname}[1]
  {\smash{\makebox[0pt]{\hspace{-2.0in}\raisebox{8pt}{\tiny #1}}}}

\thanks{The research described in this publication was made possible in part 
by Award No. RM1-2089 of the U.S. Civilian Research \& Development 
Foundation for the Independent States of the Former Soviet Union (CRDF)}

\begin{document}

\title{Coherent sheaves on configuration schemes.}

\author{Valery A.~Lunts}
\address{Department of Mathematics, Indiana University,
Bloomington, IN 47405, USA}
\email{vlunts@@indiana.edu}

\begin{abstract}
We introduce and study {\it configuration schemes}, which are obtained 
by ``glueing'' usual schemes along closed embeddings. The category 
of coherent sheaves on a configuration scheme is investigated. 
{\it Smooth} configuration schemes provide alternative 
``resolutions of singularities'' of usual singular schemes. We 
consider in detail the case when the singular scheme is a 
union of hyperplanes in a projective space. 
\end{abstract}

\maketitle

\section{Introduction}

We introduce and study {\it configuration schemes}, or rather, coherent 
sheaves on them. Namely, let $S$ be a finite poset considered as a category, 
i.e. for $\alpha ,\beta \in S$ the set of morphisms $Mor(\alpha, \beta)$ 
contains a unique element if $\alpha \geq \beta$ and is empty otherwise. 
Consider a functor $X$ from the opposite category $S^{opp}$ to the category 
of schemes. 
Thus, in particular, 
for each $\alpha \in S$ we are given a scheme $X(\alpha)$ with a morphism 
$f_{\beta \alpha }:X(\beta )\to X(\alpha)$ if $\beta \leq \alpha$. We require 
$f_{\beta \alpha}$'s to be {\it closed embeddings} and then call this data 
a configuration scheme $X/S$. (One could relax this condition by 
requiring only 
that $f_{\beta \alpha }$'s be affine morphisms). We say that $X/S$ is smooth 
if each $X(\alpha)$ is a smooth scheme.

A {\it  coherent sheaf} $F$ on $X/S$ is a collection of coherent sheaves 
$\{ F_{\alpha }\in coh(X(\alpha ))\}$ with morphisms 
$\sigma _{\alpha \beta}:f_{\beta \alpha}^*F_{\alpha}\to F_{\beta }$ 
satisfying an obvious compatibility condition. The category $coh(X/S)$ of
 coherent sheaves on $X/S$ is an abelian category. If $X/S$ is smooth, 
the category $coh(X/S)$ has finite cohomological dimension (which is usually 
twice what one expects in the coherent theory).  

Configuration schemes appear naturally in the following situation. 
Let $X_0$ be a (singular) reducible scheme with smooth components, which 
intersect along smooth subschemes. Let $S$ be the poset indexing all nonempty 
intersections of these components, and for $\alpha \in S$ let $X(\alpha )$ be
 the corresponding intersection itself. We get a smooth configuration 
scheme $X/S$ with the obvious inverse image functor (or ``localization'' 
functor)
$$\cL :coh(X_0)\to coh(X/S).$$
The configuration scheme $X/S$ can be considered as a (nontraditional) 
resolution of singularities of $X_0$, which may behave better than the 
usual resolution ${\tilde X}_0 \to X_0$. For example, one can reconstruct 
the scheme $X_0$ from the configuration scheme $X/S$, whereas $\tilde{X}_0$  
does not ``remember'' that it is a resolution of $X_0$. 
In the last section of 
this paper we give examples of some interesting auto-equivalences 
of the derived category $D^b(coh(X_0))$ and show that they ``extend'' 
to auto-equivalences of $D^b(coh(X/S))$ (whereas they do not extend to 
$D^b(coh({\tilde X}_0))$). 

The derived categories $D^b(coh(X/S))$ for smooth configuration schemes 
$X/S$ tend to have similar prorerties to those of the derived categories 
$D^b(coh(Y))$ for smooth schemes $Y$. For example, if each $X(\alpha )$ 
is a (smooth) projective variety, then every cohomological functor 
$$D^b(coh(X/S))\to Vect$$
is representable. In particular, $D^b(coh(X/))$ has a Serre functor ([BK]). 

We work out in some detail the case when $X_0$ is the union of hyperplanes 
 in general position in a projective space. We call the corresponding 
configuration scheme $X/S$ the {\it hyperplane configuraion scheme}. 
We prove that the derived localization functor $\bbL\cL:D^b(coh(X_0))\to 
D^-(coh(X/S))$ is fully faithful when restricted to the full 
subcategory of perfect complexes in $coh(X_0)$ (Theorem 4.5), which again 
shows the intimate relation of the singular scheme $X_0$ and its 
``desingularization'' $X/S$.  
Also we find a strong 
exceptional collection $\cE$ in $D^b(coh(X/S))$ ([Bo]) which implies that 
the category $D^b(coh(X/S))$ is equivalent to the derived category 
$D^b((A-mod)_f)$ of finite right $A-modules$ for a finite dimensional 
algebra $A$. This is again in analogy with $D^b(coh(Y))$ for certain 
smooth projective schemes $Y$ (see [Bo]). 

Finally I want to mention that some configuration schemes $X/S$ 
(or rather categories $coh(X/S)$) are intimately related to categories 
of {\it constructible} sheaves on algebraic varieties. In the next paper 
I will use configuration schemes to construct a mirror-type 
correspondence between coherent and constructible sheaves. 

It is my pleasure to thank the participants of the seminar in the 
Moscow Independent University, where this work was discussed. Special 
thanks are due to A.~Bondal, V.~Golyshev and D.~Orlov for their 
help and encouragement. Orlov critisized my original false conjectures 
and Bondal made some useful suggestions (including the last agrument 
in the proof of Theorem 4.5). Also Bondal independently had ideas 
similar to those of the present work.

\section{Configuration categories}

\subsection{Poset categories and configuration categories}
Let $S$ be a poset, which we consider as a category. That is for $\alpha, \beta \in 
S$ the set $Mor(\alpha, \beta)$ consists of a single element if $\alpha \geq \beta $ and 
is empty otherwise. Let $\cC$ be a functor from $S$ to the category $Cat$ of all 
categories. 
That is for each $\alpha \in S$ we have a category $\cC(\alpha)$ and if $\alpha \geq 
\beta$ 
we have a functor $\phi _{\alpha \beta}:\cC(\alpha)\to \cC(\beta)$ such that 
$\phi _{\alpha \gamma}=\phi _{\beta \gamma}\phi _{\alpha \beta}$ in case $\alpha \geq 
\beta 
\geq \gamma$. 

The pair $(S,\cC )$ defines a new category $\cM (S,\cC)$ which may be considered as a 
glueing 
of categories $\cC (\alpha )$ along the functors $\phi _{\alpha \beta}$. Namely, 
objects of $\cM (S,\cC)$ are collections 
$$F=(F_{\alpha }\in \cC(\alpha ))_{\alpha \in S}$$
together with morphisms 
$$\sigma _{\alpha \beta}(F)=\sigma _{\alpha \beta}:
\phi _{\alpha \beta}(F_{\alpha})\to F_{\beta }$$
such that $\sigma _{\alpha \gamma }=\sigma _{\beta \gamma}\cdot \phi _{\beta \gamma}
(\sigma _{\alpha \beta})$. 
Morphisms $f\in Mor(F,G)$ are collections $f=\{f_{\alpha }\in Mor_{\cC(\alpha)}(F_{\alpha},
G_{\alpha})\}$, which are compatible with the structure morphisms $\sigma _{\alpha \beta}$.

Note that each category $\cC(\alpha )$ is naturally a full subcategory of $\cM(S,\cC)$. 

\begin{defn} We call $\cM (S,\cC)$ a {\it poset category} over $S$. 
We say that $\cM (S,\cC)$ is an {\it abelian poset category} if each $\cC(\alpha )$ 
is an abelian category and the functors $\phi _{\alpha \beta}$ are 
{\it right exact}.
\end{defn}

\begin{lemma} An abelian poset category $\cM(S,\cC)$ is an abelian category.
\end{lemma}

\begin{pf} Indeed, given two objects $F, G$ in $\cM(S,\cC)$ and a morphism 
$f:F\to G$, define $Ker(f)$ and $Coker(f)$ componentwise. Namely, put 
$Ker(f)_{\alpha }:=Ker(f_{\alpha })$, $Coker(f)_{\alpha }:=Coker(f_{\alpha })$. 
Note that $Coker(f)$ is a well defined object in $\cM(S,\cC)$ since the functors 
$\phi _{\alpha \beta}$ are right exact. 
\end{pf}

In this paper we will only be interested in abelian poset categories.
For simplicity we only consider $\bbC$-linear abelian categories.

\begin{defn} An abelian poset category $\cM(S,\cC)$ is called a 
{\it configuration category} if the set $S$ is finite and 
each functor $\phi _{\alpha \beta }$ 
has a right adjoint $\psi _{\beta \alpha}$, which is exact. 
\end{defn}

\begin{example} Our main example of an abelian poset category comes from geometry. 
Namely, let $X:S^{opp}\to Sch/\bbC$ be a functor from the opposite category 
$S^{opp}$ to the category of schemes (say, of finite type over $\bbC$). 
We call 
this data a {\it poset of schemes} and denote it $X/S$. 
Then for $\alpha \geq \beta$ 
in $S$ 
we have a morphism of schemes $f_{\beta \alpha}:X(\beta )\to X(\alpha)$. 
This induces the inverse image functor $\phi _{\alpha \beta}:=f_{\beta \alpha}^*:
coh(X(\alpha ))\to coh(X(\beta )$ between the categories of coherent sheaves. 
This functor is right exact. Thus we obtain a functor $\cC:S\to Cat$, such that 
$\cC(\alpha )=coh(X(\alpha))$ and $\phi _{\alpha \beta}$ as above, which gives rise 
to the corresponding abelian poset category $\cM(S,\cC)$. An object $F$ of 
$\cM(S,\cC)$ is a collection of coherent sheaves $F=\{F_{\alpha }\in coh(X(\alpha ))\}$ 
with morphisms $\sigma _{\alpha \beta }:f^*_{\alpha \beta}F_{\alpha}\to F_{\beta }$, 
satisfying the compatibility condition. 
\end{example}

\begin{defn} In the notation of the previous example 
assume that the set $S$ is finite and all morphisms $f_{\beta \alpha}$ 
are closed embeddings. Then we call the poset of schemes $X/S$ a 
{\it configuration scheme}. 
\end{defn}

\begin{defn} Let $X/S$ be a configuration scheme. 
 The direct image functor $\psi _{\beta \alpha}:=
f_{\beta \alpha *}:coh(X(\beta))\to coh(X(\alpha ))$ is the right adjoint to $\phi 
_{\alpha \beta}$ and it is exact. Thus $\cM(S,\cC)$ becomes a configuration category, 
which is called the  {\it category of coherent sheaves on $X/S$} 
and denoted 
$coh(X/S)$. Similarly, we define the {\it category $qcoh(X/S)$ 
of quasicoherent sheaves on $X/S$}. 
Clearly, $coh(X/S)$ is a full subcategory of $qcoh(X/S)$. 
We say that $X/S$ is a {\it smooth} configuration scheme, if all schemes $X(\alpha )$ are 
smooth. 
\end{defn}

\begin{remark} Note that we could weaken the condition that 
$f_{\alpha \beta}$' are closed embeddings and require only that they be 
affine morphisms. Then the category $coh(X/S)$ is still a configuration 
category.
\end{remark}

\subsection{Motivation} Our main interest lies in the categories $coh(X/S)$ 
and $D^b(coh(X/S))$ for 
smooth configuration schemes $X/S$. Such schemes arise from degenerations of 
families of smooth schemes. 
(More precisely, this happens when the singular fiber has smooth components which 
intersect along smooth subschemes). 
We will show that 
for a smooth configuration scheme $X/S$ the category $coh(X/S)$ inherits many of the 
good properties of the categories $coh(X(\alpha ))$. For example, $coh(X/S)$ has 
finite cohomological dimension (Corollary 3.5). 
Consider the bounded derived category $D^b(coh(X/S))$. 
If each (smooth) scheme $X(\alpha )$ is projective, then every covariant or 
contravariant cohomological functor from $D^b(coh(X/S))$ to $Vect$ is 
representable (Theorem 3.14). 
In particular there exists a Serre functor in $D^b(coh(X/S))$ (Theorem 3.18). 

\subsection{Operations in configuration categories}

Let $\cM =\cM(S,\cC)$ be an abelian poset category. For $F\in \cM$ define its 
support $\supp(F)=\{\alpha \in S\vert F_{\alpha }\neq 0\}$.

Define a topology on $S$ by taking as a basis of open sets the subsets
$U_{\alpha }=\{\beta \in S\vert \beta \leq \alpha\}$. 

Note that $Z_{\alpha }=\{ \gamma \in S\vert \gamma \geq \alpha \} $ is 
a closed subset in $S$.

Let $U\subset S$ be open and $Z=S-U$ -- the complementary closed. Let $\cM_U$ 
(resp. $\cM_Z$) be the full subcategory of $\cM$ consisting of objects $F$ with 
support in $U$ (resp. in $Z$). For every object $F$ in $\cM$ there is a natural 
short exact sequence 
$$0\to F_U\to F\to F_Z\to 0,$$
where $F_U\in \cM_U$, $F_Z\in \cM_Z$. Indeed, take 
$$(F_U)_{\alpha }=\begin{cases}
F_{\alpha}, &\ \text{if $\alpha\in U$}, \\
0,          &\ \text{if $\alpha\in Z$}.
\end{cases} 
$$
$$(F_Z)_{\alpha }=\begin{cases}
F_{\alpha }, &\ \text{if $\alpha\in Z$}, \\
0,           &\ \text{if $\alpha\in U$}.
\end{cases}
$$
We may consider $U$ (resp.$Z$) as a subcategory of $S$ and restrict the functor 
$\cC$ to $U$ (resp. to $Z$). This gives abelian poset categories 
$\cM (U,\cC\vert_U)$ and $\cM(Z,\cC\vert_Z)$. 

Denote by $j:U\hookrightarrow S$ and $i:Z\hookrightarrow S$ the inclusions. 
We get the obvious restriction functors
$$j^*=j^!:\cM\to\cM(U,\cC\vert_U),\quad i^*:\cM \to \cM(Z,\cC\vert_Z).$$
Clearly these functors are exact. The functor $j^*$ has an exact left adjoint 
$j_!:\cM(U,\cC\vert_U)\to \cM$ (``extension by zero''). Its image is the 
subcategory $\cM_U$. The functor $i^*$ has an exact right adjoint 
$i_*=i_!:\cM(Z,\cC\vert_Z)\to \cM$ (also ``extension by zero''). Its image 
is the subcategory $\cM_Z$. It follows that $j^*$ and $i_*$ preserve injectives 
(as right adjoints to exact functors). We have $j^*j_!=Id$, $i^*i_*=Id$. 

Note that the short exact sequence above is just 
$$0\to j_!j^*F\to F\to i_*i^*F\to 0,$$
where the two middle arrows are the adjunction maps. 

For $\alpha \in S$ denote by $j_{\alpha }:\{\alpha \}\hookrightarrow S$ the inclusion. 
The inverse image functor $j_{\alpha }^*:\cM \to \cM _{\{\alpha \} }=
\cC(\alpha )$, $F\mapsto F_{\alpha }$ 
has a right-exact left adjoint $j_{\alpha +}$ defined as follows
$$(j_{\alpha +}P)_{\beta }=\begin{cases}
\phi _{\alpha \beta}P, &\ \text{if $\beta \leq \alpha$}, \\
0, &\ \text{otherwise}.
\end{cases}
$$
Thus for $P\in \cC(\alpha )$, $\supp j_{\alpha +}P\subset U_{\alpha }$. 

We also consider the ``extension by zero'' functor 
$j_{\alpha !}:\cC (\alpha )\to \cM $ defined by 
$$j_{\alpha !}(P)_{\beta}=\begin{cases}
P, &\ \text{if $\alpha =\beta$}, \\
0, &\ \text{otherwise}.
\end{cases}
$$

\begin{lemma} Assume that $\cM(S,\cC)$ is a configuration category. 
Then the  functor $j_{\alpha }^*:\cM (S,\cC)\to \cC (\alpha )$ has a 
right adjoint $j_{\alpha *}$. The functor $j_{\alpha *}$ is exact 
and preserves injectives. For $P\in \cC(\alpha )$ $\supp (j_{\alpha *}P)
\subset Z_{\alpha }$.
\end{lemma}

\begin{pf}
Given $P\in \cC (\alpha )$ we set 
$$j_{\alpha *}(P)_{\gamma}=\begin{cases}
\psi_{\alpha \gamma}(P), &\ \text{if $\gamma \geq \alpha$}, \\
0,                      &\ \text{otherwise},
\end{cases}
$$
and the structure map
$$\sigma _{\gamma \delta } :\phi _{\gamma \delta }(j_{\alpha *}P)_{\gamma }
\to (j_{\alpha *}P)_{\delta }$$
is the adjunction map
$$\phi _{\gamma \delta }\psi _{\alpha \gamma }P=
\phi _{\gamma \delta }\psi _{\delta \gamma }\psi _{\alpha \delta }P
\to \psi _{\alpha \delta }P$$
if $\alpha \leq \delta \leq \gamma$ and $\sigma _{\gamma \delta }=0$ 
otherwise. 

It is clear that $j_{\alpha *}$ is exact (since 
$\psi $'s are such) and that $\supp (j_{\alpha *}P)\subset Z_{\alpha }$. 

Let us prove that $j_{\alpha *}$ is the right adjoint to $j_{\alpha}^*$. 

Let $P\in \cC(\alpha )$ and $M\in \cM (S,\cC )$. Given $f_{\alpha }\in 
\Hom (M_{\alpha },P)$ for each $\gamma \geq \alpha $ we obtain a map 
$f_{\alpha}\cdot\sigma _{\gamma \alpha }(M):
\phi _{\gamma \alpha }M_{\gamma }\to P$ 
and hence by adjunction $f_{\gamma }:M_{\gamma }\to \psi _{\alpha \gamma }P=
(j_{\alpha *}P)_{\gamma }$. The collection $f=\{f_{\gamma }\} $ is a morphism 
$f:M\to j_{\alpha *}P$. It remains to show that the constructed map 
$$\Hom (M_{\alpha },P)\to \Hom (M, j_{\alpha *}P)$$
is surjective or, equivalently, that the restriction map 
$$\Hom (M,j_{\alpha *}P)\to \Hom (M_{\alpha },P),\quad f\mapsto f_{\alpha }$$
is injective. 

Assume that  $0\neq f\in \Hom (M,j_{\alpha *}P)$, i.e. $f_{\gamma }\neq 0$ 
for some $\gamma \geq \alpha $. By definition we have the commutative 
diagram 
$$
\begin{CD}
\phi _{\gamma \alpha }M_{\gamma }@ >\phi_{\gamma \alpha}(f_{\gamma })>> 
\phi _{\gamma \alpha }\psi _{\alpha \gamma }P\\
@V\sigma _{\gamma \alpha }(M)VV  @VV\epsilon_PV\\
M_{\alpha }@>f_{\alpha }>>P,
\end{CD}
$$
where $\epsilon_P$ is the adjunction morphism. Note that 
$\epsilon_P\phi _{\gamma \alpha }(f_{\gamma }):\phi _{\gamma \alpha}
M_{\gamma }\to P$ is the morphism, which corresponds to 
$f_{\gamma }: M_{\gamma }\to \psi _{\alpha \gamma }P$ by the adjunction 
property. Hence $\epsilon _P\phi _{\gamma \alpha }(f_{\gamma })\neq 0$. 
Therefore $f_{\alpha }\neq 0$. This shows the injectivity of the restriction 
map $f\mapsto f_{\alpha }$ and proves that $j_{\alpha *}$ is the right adjoint 
to $j_{\alpha }^*$. Finally, $j_{\alpha *}$ preserves injectives being the 
right adjoint to an exact functor. 
\end{pf}

\subsection{Cohomological dimension of configuration categories}
\begin{lemma}
Let $\cM=\cM(S,\cC)$ be a configuration category. Assume that each category 
$\cC(\alpha )$ has enough injectives. Then so does $\cM$. 
\end{lemma}

\begin{pf}
We will prove the assertion by induction on the cardinality of $S$. 
Let $\beta \in S$ be a smallest element. Put $U=U_{\beta }=\{\beta \}$, 
$Z=S-U$. Let $j=j_{\beta}:U\hookrightarrow S$ and $i:
Z\hookrightarrow S$ be the corresponding 
open and closed embedding. 

For $F$ in $\cM$ we need to find a monomorphism $F\to I$, where $I$ is an 
injective object in $\cM$. It suffices to do so for $F_U=j_!j^*F$ and $F_Z=i_*i^*F$. 
Let $j^*F\to I_1$, $i^*F\to I_2$ be similar monomorphisms in categories 
$\cM(U,\cC\vert_U)=\cC(U)$ and $\cM(Z,\cC\vert_Z)$ 
(this is possible by the induction 
assumption). Then 
$$j_*j^*F\to j_*I_1,\quad i_*i^*F\to i_*I_2$$
are monomorphisms in $\cM$ with $j_*I_1$ and $i_*I_1$ being injective. It remains to 
note the obvious monomorphism $j_!j^*F\to j_*j^*F$. 
\end{pf}

\begin{prop}
Let $\cM=\cM(S,\cC)$ be a configuration category. Assume that each category 
$\cC(\alpha )$ has enough injectives and has finite cohomological dimension. Then 
$\cM$ also has finite cohomological dimension. 
\end{prop}

\begin{pf}
We argue by induction on the cardinality of $S$ as in the proof of the last lemma. 
In fact we will use the notation of that proof. 

Given $F$ in $\cM$ we need to find a finite injective resolution of $F$. It suffices 
to do so for $j_!j^*F$ and $i_*i^*F$. Let $j^*F\to I_1^{\bullet}$, $i^*F\to 
I_2^{\bullet}$ be such resolutions in categories $\cM(U,\cC\vert_U)$ and 
$\cM(Z,\cC\vert_Z)$ respectively. Then $j_*j^*F\to j_*I_1^{\bullet}$, 
$i_*i^*F\to i_*I_2^{\bullet}$ will be injective resolutions in $\cM$. Consider the 
short exact sequence 
$$0\to j_!j^*F\to j_*j^*F\to G\to 0.$$
Then $\supp(G)\subset Z$ and so by induction $G=i_*i^*G$ has a finite injective resolution 
in $\cM$. Therefore the same is true for $j_!j^*F$.
\end{pf}

\subsection{Derived categories of configuration categories}
Let $\cM=\cM(S,\cC)$ be a configuration category. Denote by $D^b(\cM)$ its bounded derived 
category.

Let $U\stackrel{j}{\hookrightarrow}S\stackrel{i}{\hookleftarrow}Z$ be embeddings of an 
open $U$ and a complementary closed $Z$. The exact functors $j^*,j_!,i^*,i_*$ extend 
trivially to corresponding functors between derived categories $D^b(\cM)$, 
$D^b(\cM(U,\cC\vert_U))$, $D^b(\cM(Z,\cC\vert_Z))$. 

If $j_{\alpha }:\{\alpha \}\hookrightarrow S$ is the inclusion, 
then the same is true for the exact functor $j_{\alpha *}$. When this 
causes no confusion we will denote 
the derived functor again by the same symbol, i.e. we will write $j_!$ for 
$\bbR j_!$, etc. 

Recall that functors $j_!$ and $i_*$ identify categories $\cM(U,\cC\vert_U)$ and 
$\cM(Z,\cC\vert_Z)$ with $\cM_U$ and $\cM_Z$ respectively. Denote by 
$D_U^b(\cM)$ and $D^b_Z(\cM)$ the full subcategories of $D^b(\cM)$ consisting of complexes 
with cohomology objects in $\cM_U$ and $\cM_Z$ respectively. 

\begin{lemma} Assume that all categories $\cC(\alpha )$ have enough injectives. Then 
in the above notation the functors $i_*:D^b(\cM(Z,\cC\vert_Z))\to D^b(\cM)$ and 
$j_!:D^b(\cM(U,\cC\vert_U))\to D^b(\cM)$ are fully faithful. The essential images of these 
functors are the full subcategories $D_Z^b(\cM)$ and $D_U^b(\cM)$ respectively. 
\end{lemma}

\begin{pf}
We need to prove that functors $i_*$ and $j_!$ preserve $\Ext$'s between objects. 

Let $F,G\in \cM(Z,\cC\vert_Z)$ and let $G\to I^{\bullet}$ be an injective resolution 
(it exists by Lemma 2.9). Then $i_*G\to i_*I^{\bullet}$ is still an injective resolution. 
We have 
\begin{eqnarray*}
\Ext^k(i_*F,i_*G) & = & H^k\Hom^{\bullet}(i_*F,i_*I^{\bullet}) \\
                  & = & H^k\Hom^{\bullet}(i^*i_*F,I^{\bullet}) \\
                  & = & H^k\Hom^{\bullet}(F,I^{\bullet}) \\
                  & = & \Ext^k(F,G).
\end{eqnarray*}

Thus $i_*$ preserves $\Ext$'s. 
  
Let $F,G\in \cM(U,\cC\vert_U)$. Let $j_!G\to J^{\bullet}$ be an injective resolution. 
Then $G=j^*j_!G\to j^*J^{\bullet}$ is still an injective resolution. We have 
\begin{eqnarray*}
\Ext^k(j_!F,j_!G) & = & H^k\Hom^{\bullet}(j_!F,J^{\bullet}) \\
                  & = & H^k\Hom^{\bullet}(F,j^*J^{\bullet}) \\
                  & = & \Ext^k(F,G).
\end{eqnarray*}

\end{pf}

We immediately obtain the following corollary

\begin{cor} 
Assume that all categories $\cC(\alpha )$ have enough injectives. Then the categories 
$D^b(\cM_U)$ and $D^b(\cM_Z)$ are naturally equivalent to $D^b_U(\cM)$ and 
$D_Z^b(\cM)$ respectively.
\end{cor}

\begin{cor}
Assume that each category $\cC(\beta )$ has enough injectives. Fix 
$\alpha \in S$. Let $i:\{\alpha \}\hookrightarrow U_{\alpha }$ and 
$j:U_{\alpha }\hookrightarrow S$ be the closed and the open embeddings 
respectively. Then the composition 
functor 
$$j_!i_*:D^b(\cC(\alpha ))\to D^b(\cM (S,\cC))$$
is fully faithful. In particular, the derived category $D^b(\cC(\alpha ))$ 
is naturally a full subcategory of $D^b(\cM (S,\cC))$. 
\end{cor}

\begin{pf}
Indeed, by Lemma 2.11 above the functors 
$$i_*:D^b(\cC(\alpha ))\to D^b(\cM (U_{\alpha },\cC\vert _{U_{\alpha }}))$$ 
and 
$$j_!:D^b(\cM (U_{\alpha },\cC\vert _{U_{\alpha }}))\to D^b(\cM (S,\cC))$$ 
are fully faithful. So is their composition.
\end{pf}

Recall the following definitions from [BK]. 

\begin{defn}
Let $\cA$ be an additive category, $\cB\subset \cA$ -- a full subcategory. A 
{\it right orthogonal} to $\cB$ in $\cA$ is a full subcategory $\cB^{\perp}\subset \cA$ 
consisting of all objects $C$ such that $\Hom(B,C)=0$ for all $B\in \cB$.
\end{defn}

\begin{defn}
Let $\cA$ be a triangulated category, $\cB \subset \cA$ -- a full triangulated subcategory. 
We say that $\cB$ is {\it right-admissible} if for each $X\in \cA$ there exists an 
exact triangle $B\to X\to C$ with $B\in \cB$, $C\in \cB^{\perp}$.
\end{defn}

\begin{lemma}
a) In the category $\cM$ we have $\cM_U^{\perp}=\cM_Z$.

b) Assume that all categories $\cC(\alpha )$ have enough injectives (so that, 
in particular, $D^b(\cM_U)$ and $D^b(\cM_Z)$ are full subcategories of $D^b(\cM)$). 
Then

i) in $D^b(\cM)$ we have $D^b(\cM_U)^{\perp}=D^b(\cM_Z)$,

ii) the subcategory $D^b(\cM_U)\subset D^b(\cM)$ is right-admissible.

\end{lemma}

\begin{pf}
a). Given $F\in \cM_U$, $G\in \cM$ we have $\Hom(F,G)=\Hom(F,G_U)$. Hence 
$\Hom(F,G)=0$ for all $F$ iff $G_U=0$ or, equivalently, $G\in \cM_Z$. 

b)i). Let $G^{\bullet}\in D^b(\cM)$. If $G^{\bullet}\in D^b(\cM_U)^{\perp}$, 
then 
 $G^{\bullet}_U$ is acyclic, i.e. $G^{\bullet}\in D^b_Z(\cM)$. Vice versa, 
an injective complex $G^{\bullet}$ in $D^b(\cM_Z)$ stays injective in $D^b(\cM)$. 
Thus $D^b(\cM_Z)\subset D^b(\cM_U)^{\perp}$.

ii). Given $X^\bullet \in D^b(\cM)$ the required exact triangle is  
$X^\bullet _U\to X^\bullet \to X^\bullet _Z$. 
\end{pf}

\section{Categories $coh(X/S)$ and $qcoh(X/S)$ for smooth configuration schemes $X/S$}

We now turn to the subject of our primary interest: the category $coh(X/S)$ for a smooth 
configuration scheme $X/S$. This category does not have enough injectives, so we will 
study $coh(X/S)$ (as usual) by embedding it in a larger category $qcoh(X/S)$.

So let $X/S$ be a smooth configuration scheme.

\begin{lemma} The category $qcoh(X/S)$ has enough injectives.
\end{lemma}

\begin{pf} This is a direct consequence of the general Lemma 2.9.
\end{pf}

\begin{cor} The category $qcoh(X/S)$ has finite cohomological dimension.
\end{cor}

\begin{pf} This follows from Lemma 3.1, Proposition 2.10 and the well 
known fact that the category $qcoh(Y)$ has finite cohomological dimension 
if $Y$ is smooth.
\end{pf}

Consider the bounded derived categories $D^b(coh(X/S))$ and $D^b(qcoh(X/S))$. 
Let $D^b_{coh}(qcoh(X/S))\subset D^b(qcoh(X/S))$ be the full subcategory 
consisting of complexes with cohomologies in $coh(X/S)$. The next proposition 
is the analogue of a similar well known result for usual schemes. 

\begin{prop}
The categories $D^b(coh(X/S))$ and $D^b_{coh}(qcoh(X/S))$ are equivalent.
\end{prop}

\begin{lemma}
Let $f:F\to T$ be an epimorphism, where $F\in qcoh(X/S)$ and $T\in coh(X/S)$. 
There exists a subobject $G\subset F$ such that $G\in coh(X/S)$ and 
$f\vert_G:G\to T$ 
is still an epimorphism.
\end{lemma}

\begin{pf}
The lemma is well known for usual noetherian schemes. For each $\alpha \in S$ 
choose a coherent subsheaf $P_{\alpha }\subset F_{\alpha }$ such that 
$f_{\alpha}\vert_{P_{\alpha}}:P_{\alpha}\to T_{\alpha}$ is an epimorphism. 
The inclusion $P_{\alpha }\hookrightarrow F_{\alpha}$ defines the corresponding morphism 
$t_{\alpha}:i_{\alpha +}P_{\alpha }\to F$. Put 
$$P:=\bigoplus_{\alpha}i_{\alpha +}P_{\alpha}\in coh(X/S),\quad 
t=\oplus t_{\alpha} :P\to F.$$
Then the composition $f\cdot t:P\to T$ is an epimorphism and we can take $G=Im(t)$. 
\end{pf}

\begin{pf}(of proposition). The proposition follows from the above lemma by 
standard homological algebra as (for example) in [ST] (p.12,13).
\end{pf}

\begin{cor} The category $coh(X/S)$ has finite cohomological dimension.
\end{cor}

\begin{pf} This follows from Corollary 3.2 and Proposition 3.3.
\end{pf}

Let $U\stackrel{j}{\hookrightarrow}S\stackrel{i}{\hookleftarrow}Z$ be embeddings 
of an open and a complementary closed subsets. Consider the functors 
$$j_!:D^b(coh(X/U))\to D^b(coh(X/S)),$$ 
$$i_*:D^b(coh(X/Z))\to D^b(coh(X/S)).$$

\begin{prop}
i) The functors $j_!$ and $i_*$ are fully faithful. Thus $D^b(coh(X/U))$ and 
$D^b(coh(X/Z))$ can be considered as full 
triangulated subcategories of $D^b(coh(X/S))$. 

ii) $D^b(coh(X/Z))$ is the right orthogonal to $D^b(coh(X/U))$ in $D^b(coh(X/S))$. 

iii) $D^b(coh(X/U))$ is right-admissible in $D^b(coh(X/S))$. 

\end{prop}

\begin{pf} i). By Lemma 2.11 the analogous fact holds for 
the categories $D^b(qcoh)$ instead of 
$D^b(coh)$. Hence it also holds for the categories $D^b_{coh}(qcoh)$. Therefore it 
holds for $D^b(coh)$ by the above proposition. 

ii). By Lemma 2.16 b)i) the analogous fact holds for the 
categories $D^b(qcoh)$ instead of $D^b(coh)$. 
Actually, the same proof shows that the assertion is true for categories 
$D^b_{coh}(qcoh)$. Hence by Proposition 3.3 it is also true for $D^b(coh)$. 

iii). The required exact triangle is $F_U\to F\to F_Z$. 
\end{pf}

\begin{cor}
Fix $\alpha \in S$. Let  $i:\{\alpha \}\hookrightarrow U_{\alpha}$ and 
$j:U_{\alpha }\hookrightarrow  S$ be closed and the open embeddings 
respectively. Then the composition of functors 
$$j_!i_*:D^b(coh(X(\alpha )))\to D^b(coh(X/S))$$
is fully faithful. In particular, the derived category $D^b(coh(X(\alpha )))$ 
is naturally a full subcategory of $D^b(coh(X/S))$.
\end{cor}

\begin{pf} By Proposition 3.6 above the functors $j_!$ and $i_*$ are fully 
faithful. So is their composition.
\end{pf}

Recall the following definitions from [BK].

\begin{defn} A triangulated category $\cA$ is of {\it finite type} if for all 
$B,C\in \cA$ the space $\Hom(B,C)$ is finite dimensional and 
$\Ext ^i(B,C)=0$ for $\vert i \vert >>0$. 
\end{defn}

\begin{defn} Let $\cA$ be a triangulated category and $h:\cA \to Vect$ be a contravariant 
cohomological functor. We say that $h$ is of {\it finite type} if for all $E\in \cA$ 
the vector 
spaces $h^i(E)=h(E[-i])$ are finite dimensional and are equal to zero for almost all $i$. 
Similarly, for covariant cohomological functors. 
\end{defn}

\begin{defn} Let $\cA$ be a triangulated category of finite type. We say that $\cA$ is 
{\it right-saturated} if every contravariant cohomological functor 
$h:\cA \to Vect$ of finite type is representable, i.e. $h(\cdot)=\Hom(\cdot ,X)$ 
for some $X\in \cA$. Similarly we define left-saturatedness using covariant 
cohomological functors. The category $\cA$ is {\it saturated} if it is left- and 
right-saturated. 
\end{defn}

One of the main results in [BK] is the following theorem.

\begin{thm} Let $X$ be a smooth projective variety. Then the category $D^b(coh(X))$ 
is saturated. 
\end{thm}

The technique developed in [BK] allows us to extend the above theorem to the category 
$D^b(coh(X/S))$ for a smooth configuration scheme with projective $X(\alpha )$'s. 

\begin{lemma} Assume that each $X(\alpha )$ is a (smooth) projective variety. 
Then the 
category $D^b(coh(X/S))$ is of finite type.
\end{lemma} 

\begin{pf} Put $\cM=coh(X/S)$. Let $\beta \in S$ be a smallest element, $U=U_{\beta }=
\{\beta \}$ -- the corresponding open subset, $Z=S-U$ -- the complementary closed. 
Given $F,G\in \cM$ consider the short exact sequences in $\cM$ 
$$0\to F_U\to F\to F_Z\to 0,\quad 0\to G_U\to G\to G_Z\to 0.$$
To show that $\oplus_kdim\Ext^k(F,G)< \infty$ it suffices to prove the same for the 
other four members of the above exact sequences. 
By induction on the cardinality of $S$ 
(and using Proposition 3.6,i) we may assume that 
$$\bigoplus_kdim\Ext^k(F_U,G_U)<\infty,\quad \bigoplus_kdim\Ext^k(F_Z,G_Z)<\infty.$$
Also by Proposition 3.6,ii, $\Ext^{\bullet}(F_U,G_Z)=0$. It remains to consider 
$\Ext^k(F_Z,G_U)$. 

Denote by $j:U\hookrightarrow S$ the open embedding. Consider the short 
exact sequence 
$$0\to G_U\to j_*j^*G_U\to P\to 0.$$
Then $P\in \cM_Z$ and by induction $\oplus_kdim\Ext^k(F_Z,P)<\infty$. 
We claim that $\Ext^{\bullet}(F_Z,j_*j^*G_U)=0$. Indeed, let $j^*G_U\to I^{\bullet}$ be 
an injective resolution (in $qcoh(X(\beta ))$). Then $j_*j^*G_U\to j_*I^{\bullet}$ remains 
an injective resolution (Lemma 2.8), so that $\Ext^{\bullet}(F_Z,j_*j^*G_U)=
H^{\bullet}\Hom^{\bullet}(F_Z,j_*I^{\bullet})$. But for any 
$T\in qcoh(X(\beta))$  
$$\Hom(F_Z,j_*T)=\Hom(j^*F_Z,T)=\Hom(0,T)=0.$$ 
Hence also $\oplus_kdim\Ext^k(F_Z,G_U)<\infty$. 
\end{pf}

The next theorem is again from [BK].

\begin{thm} Let $\cA$ be a triangulated category of finite type, 
$\cB\subset \cA$ - a right-admissible triangulated subcategory, $\cC=\cB^{\perp}$. 
Assume that $\cB$ and $\cC$ are right- (resp. left-) saturated. Then so is $\cA$. 
\end{thm}

We have the following counterpart of Theorem 3.11.

\begin{thm} Assume that each $X(\alpha )$ is a (smooth) projective variety. Then the 
category $D^b(coh(X/S))$ is saturated.
\end{thm}

\begin{pf}  We will use Theorem 3.11 as a basis for induction on the cardinality  of $S$. 
By Lemma 3.12 the category $D^b(coh(X/S))$ is of finite type. Let $U\subset S$ be open and 
$Z=S-U$. By Proposition 3.6 the full subcategory $D^b(coh(X/U))\subset D^b(coh(X/S))$ is right-admissible 
and $D^b(coh(X/U))^{\perp}=D^b(coh(X/Z))$. By induction hypothesis the categories 
$D^b(coh(X/U))$ and $D^b(coh(X/Z))$ are saturated. Hence by Theorem 3.13 so is $D^b(coh(X/S))$. 
\end{pf}

\subsection{The Serre functor}

\begin{defn} ([BK]) Let $\cA$ be a triangulated category with finite dimensional 
$\Hom$'s. Let $F:\cA\to \cA$ be a covariant auto-equivalence which commutes with 
the translation functor [1]. We call $F$ a {\it Serre functor} if there exist 
bifunctorial isomorphisms 
$$\phi _{E,G}:\Hom (E,G)\to \Hom(G,F(E))^*$$
(for $E,G\in \cA$) with the property that the 
composition 
$$(\phi ^{-1}_{F(E),F(G)})^*\cdot \phi _{E,G}:\Hom (E,G)\to \Hom(G,F(E))^*\to 
\Hom(F(E),F(G))$$
coincides with the isomorphism induced by $F$.
\end{defn}

\begin{example} Let $X$ be a smooth projective variety, $n=dimX$, $\cA=D^b(coh(X))$, 
$K_X=\Omega _X^n$ - the canonical sheaf. Then $F(\cdot)=
(\cdot)\otimes K_X[n]$ is 
a Serre functor in $\cA$.
\end{example}

\begin{prop} ([BK]) Let $\cA$ be a triangulated category of finite type.

i) A Serre functor $F:\cA \to \cA$ is exact, i.e. it takes exact triangles to 
exact triangles.

ii) Any two Serre functors in $\cA$ are canonically isomorphic.

iii) If $\cA$ is saturated, then it has a Serre functor.
\end{prop}

We immediately obtain the following theorem.

\begin{thm} Assume that each $X(\alpha )$ is a (smooth) projective variety. 
Then the category $D^b(coh(X/S))$ has a Serre functor.
\end{thm}

\begin{pf} Follows from Theorem 3.14 and Proposition 3.17.
\end{pf}

\subsection{Cohomology of coherent sheaves on configuration schemes}

\medskip

Let $X/S$ be a configuration scheme.

\begin{defn} Consider the object $G\in coh(X/S)$ such that $G_{\alpha }=
\cO _{X(\alpha )}$ and for $\beta \leq \alpha$ the map 
$\sigma _{\alpha \beta }:\phi _{\alpha \beta}(G_{\alpha })\to G_{\beta }$ 
is the identity. We call this object the {\it structure sheaf of} $X/S$ 
and denote it $\cO _{X/S}$. 
\end{defn}

\begin{defn} Let $F\in coh(X/S)$. We define the cohomology  
of $F$ as 
$$H^i(F)=H^i(X/S,F):=\Ext ^i(\cO _{X/S},F).$$
\end{defn}

Note that any $G\in coh(X/S)$ has a finite left resolution by sheaves 
which are finite direct sums of $G_{U_{\alpha }}$'s for various $\alpha 
\in S$:
$$0\to \oplus G_{U_{\beta }}\to ...\to \oplus G_{U_{\alpha }}\to G\to 0.$$
Hence, in particular, $\cO _{X/S}$ has such a resolution. Thus to 
compute $H^i(X/S,F)$ we can use a spectral sequence with the $E_1$ term  
equal to a direct sum of vector spaces $H^i(X/U_{\alpha },F)$ for various 
$\alpha \in S$, $i\geq 0$. 

The next lemma shows how to compute the cohomology 
$H^\bullet(X/U_{\alpha },F)$. 

\begin{lemma}
Assume that $\alpha \in S$ is the unique maximal element, i.e. 
$S=U_{\alpha }$. Then for any $F\in coh(X/S)$ there is a natural isomorphism 
$$H^\bullet(X/S,F)\simeq H^\bullet (X(\alpha ),F_{\alpha }).$$
\end{lemma}

\begin{pf} Let $j_{\alpha }:\{\alpha \}\hookrightarrow S$ be the 
embedding. Note that $\cO _{X/S}=j_{\alpha +}\cO _{X(\alpha )}$. 
Hence for 
$F\in coh(X/S)$
$$\Hom (\cO _{X/S},F)\simeq \Hom (\cO _{X(\alpha )},F_{\alpha }),$$
which proves the lemma for $H^0$. To complete the proof of the lemma 
it suffices to prove the following claim. 

\bigskip
\noindent{\it Claim.} {\it Every $F\in coh(X/S)$ can be embedded in an 
injective  object $I\in qcoh(X/S)$ such that $I_{\alpha }$ is acyclic for 
the functor $\Hom (\cO _{X(\alpha )},\cdot)$. }

\medskip
\noindent{\it Proof of the claim.} For $\beta \in S$ let $j_{\beta }:
\{ \beta \}\hookrightarrow S$ be the inclusion. For each $\beta $ 
choose a monomorphism  $F_{\beta }\to I_{\beta }$, where 
$I_{\beta }\in qcoh(X(\beta ))$ is an injective object. By adjunction 
this defines a map
$$t^{\beta }:F\to j_{\beta *}I_{\beta },$$
where the object $j_{\beta *}I_{\beta }\in qcoh(X/S)$ 
is injective (Lemma 2.8). 
Put 
$$I:=\oplus _{\beta \in S}j_{\beta *}I_{\beta }$$
with the embedding
$$t=\oplus t^{\beta }:F\hookrightarrow I.$$

We claim that for each $\beta $ the sheaf $(j_{\beta *}I_{\beta })_{\alpha }
\in qcoh(X(\alpha ))$ is acyclic for $\Hom (\cO _{X(\alpha )},\cdot )$ (and 
hence $I$ is such). Indeed, it is equal to $\psi _{\beta \alpha }(I_{\beta })
=f_{\beta \alpha *}(I_{\beta })$ for the closed embedding 
$f_{\beta \alpha}:X(\beta )\hookrightarrow X(\alpha )$, and 
$H^\bullet(X(\beta ),I_{\beta })=H^\bullet (X(\alpha ),f_{\beta \alpha *}
I_{\beta })$. This proves the claim and the lemma.
\end{pf}

We immediately obtain the following corollary.

\begin{cor}
Assume that $\alpha \in S$ is the unique maximal element, i.e. 
$S=U_{\alpha }$. Let $F\in coh(X/S)$ be such that $F_{\alpha }=0$. Then 
$H^\bullet(X/S,F)=0$.
\end{cor}

\section{Hyperplane configuration schemes}

\subsection{Definition and cohomological dimension of hyperplane 
configuration schemes}

We want to consider in detail the simpliest type of smooth configuration 
schemes, namely, the ones defined by $d$ hyperplanes in general position 
in $\bbP^{n+1}$. 

Let $L_1,L_2,...,L_d$ be hyperplanes in general position in $\bbP^{n+1}$. 
Let $S$ be the poset of nonempty subsets of $\{1,...,d\}$ of cardinality 
$\leq n+1$. (For $\alpha$, $\beta\in S$ we put $\beta \leq \alpha$ iff 
$\beta \subset \alpha$). For $\alpha = \{i_1,...i_{j+1}\}\in S$ define 
$$X(\alpha ):=L_{i_1}\cap ...\cap L_{i_{j+1}}\simeq \bbP^{n-j}.$$
Put $d(\alpha )=dimX(\alpha )=n-j$. We will also denote 
$\bbP _{\alpha}=X(\alpha )\simeq \bbP ^{d(\alpha )}$.
If $\alpha \geq \beta $ we have the natural closed embedding 
$f_{\beta \alpha}:\bbP _{\beta }\to \bbP _{\alpha }.$ This defines a functor 
$X:S^{opp}\to Sch/\bbC$, which is a smooth configuration scheme.

\begin{defn}
We call a smooth configuration scheme $X/S$ as above {\it a hyperplane 
configuration scheme (of type $(n,d)$)}.
\end{defn}

By the general Corollary 3.5 we know that the cohomological 
dimension of the category $coh(X/S)$ is finite. The next proposition 
is a more precise result.

\begin{prop}
Let $X/S$ be a hyperplane configuration scheme of type $(n,d)$. 
Then the cohomological dimension of $coh(X/S)$ is $2n$ if $d\geq n$ and 
$n+d$ if $d<n$. I.e. it is equal to $min\{2n,\ n+d\}$.
\end{prop}

\begin{pf}
First we prove the upper bound on the cohomological dimension.

For each $k\geq 0$ denote by $Z_k\subset S$ the closed subset
$$Z_k:=\{\alpha \in S\vert d(\alpha )\geq k\}.$$
Thus $Z_n\subset Z_{n-1}\subset ...$. 

Let $F,G\in coh(X/S)$. 
Assume that 
$\supp (G)\subset Z_k$. It suffices to prove that $\Ext ^i(F,G)=0$ for 
$i>2n-k$. We will prove it by descending induction on $k$. 

If $k=n$, then 
$$\Ext^i(F,G)=\oplus_j\Ext^i(F\vert _{L_j},G\vert _{L_j}).$$
Thus $\Ext ^i(F,G)=0$ for $i>n$. 
Assume the assertion is proved for $k$ and $\supp G\subset Z_{k-1}$. 
Put $U=Z_{k-1}-Z_k$. Consider the short exact sequence 
$$0\to G_U\to G\to G_{Z_k}\to 0.$$
It suffices to prove that $\Ext ^i(F,G_U)=0$ for $i>2n-k+1$. 
We may assume that $S=Z_{k-1}$.
Note that 
$$G_U=\oplus_{\alpha \in U}G_{\alpha}.$$
Fix one $\alpha \in U$. This is a smallest element in $S$. 
Let $j_{\alpha }:\{\alpha \}\hookrightarrow S$ denote the 
open embedding. We have the short exact sequence of sheaves on $S$
$$0\to G_{\alpha }\to j_{\alpha *}j_{\alpha}^*G_{\alpha }\to 
(j_{\alpha *}j_{\alpha}^*G_{\alpha })/
G_{\alpha }\to 0.$$ 
Since $j_{\alpha *}$ is exact and preserves injectives, the object 
$j_{\alpha *}j_{\alpha }^*G_{\alpha }$ has an injective resolution of 
length $\leq k-1$. On the other hand $\supp ((j_{\alpha *}j_{\alpha}^*
G_{\alpha })/G_{\alpha })\subset Z_k$. Thus by induction 
$$\Ext ^i(F,(j_{\alpha *}j_{\alpha }^*G_{\alpha })/G_{\alpha })=0$$ 
for $i>2n-k$. Hence $\Ext^i(F,G_{\alpha })=0$ for $i>2n-k+1$. This 
completes the induction step and proves the upper bound on the 
cohomological dimension.  

To show that the bound is actually achieved choose an element 
$\beta \in S$. Put $m:=dim\bbP _{\beta }$. Choose a point 
$p\in \bbP _{\beta }$. Since $\bbP _{\beta}\subset \bbP _{\gamma}$ for 
$\beta \leq \gamma$, the point $p$ also belongs to all $\bbP _{\gamma}$ 
for $\gamma \geq \beta$. We denote the corresponding sky-scraper sheaf 
$\bbC _p^{\gamma}\in coh(\bbP _{\gamma})$. For $\gamma \in S$ denote 
the embedding $j_{\gamma }:\{\gamma \}\hookrightarrow S$. Let 
$\alpha \geq \beta$ be a maximal element. It suffices to prove that 
$$\Ext^{2n-m}(j_{\alpha *}\bbC _p^{\alpha },j_{\beta !}\bbC _p^{\beta})
\neq 0.$$
In case $m=n$, (i.e. $\alpha =\beta $) we have 
$$\Ext ^n(j_{\alpha *}\bbC _p^{\alpha },j_{\alpha !}\bbC _p^{\alpha })=
\Ext ^n_{coh(\bbP _{\alpha })}(\bbC _p,\bbC _p)\neq 0.$$

Now we argue by descending induction on $m$. The sheaf 
$j_{\beta !}\bbC _p^{\beta }$ has a resolution 
$$0\to j_{\beta !}\bbC _p^{\beta }\to j_{\beta *}\bbC _p^{\beta }\to 
\bigoplus_{\stackrel{\gamma \geq \beta}{d(\gamma )=d(\beta)+1}}
j_{\gamma *}\bbC _p^{\gamma }\to ...\to 
\bigoplus_{\stackrel{\delta \geq \beta}{d(\delta)=n}}j_{\delta *}
\bbC _p^{\delta} \to 0.$$
So it suffices to prove that $\Ext ^\bullet(j_{\alpha *}\bbC _p^{\alpha}, 
j_{\gamma *}\bbC _p^{\gamma })=0$ for all $\gamma <\alpha$. But this follows 
because $j_{\gamma *}$ is exact and preserves injectives and 
$\Hom (j_{\alpha *}\bbC _p^{\alpha},j_{\gamma *}T)=0$ for all $T\in qcoh(
\bbP _{\gamma })$. This proves the proposition.
\end{pf}

\subsection{The \v{C}ech complex}

Let $X/S$ be a hyperplane configuration scheme of type $(n,d)$. Consider 
the open covering of $S$ by the (irreducible) sets $U_{\alpha }$, $d(\alpha )
=n$. Given $F \in D^b(coh(X/S))$ we have the corresponding 
\v{C}ech resolution ${\check C}^\bullet (F)\to F$ with 
$${\check C}^\bullet(F):=0\to \bigoplus_{d(\gamma )=0}F
_{U_{\gamma }}\to ...\bigoplus_{d(\beta )=n-1}F_{U_{\beta}}\to 
\bigoplus_{d(\alpha )=n}F_{U_{\alpha}}\to 0,$$
where $F_{U_{\delta}}:=j_{\delta !}j^*_{\delta }F$ for the open 
embedding $j_{\delta }:U_{\delta }\hookrightarrow S$. The differentials 
in this complex are the sums of the natural embeddings $F_{U_{\beta}}
\hookrightarrow F_{U_{\alpha}}$ taken with $\pm$ sign. Notice that 
for each $\alpha \in S$ there is exactly one summand $F_{U_{\alpha}}$ 
in ${\check C}^\bullet(F)$. 

In case $F=\cO _{X/S}$ we will denote ${\check C}^\bullet=
{\check C}^\bullet(\cO _{X/S})$. Note that for each $\alpha \in S$, 
$j\geq 0$, $i>0$, 
$$\Ext ^i(\cO _{U_{\alpha}},\cO _{X/S}(j))=
\Ext^i(\cO _{\bbP _{\alpha }},\cO _{\bbP _{\alpha }}(j))=0.$$ 
Hence the complex $\Hom ^\bullet({\check C}^\bullet,\cO _{X/S}(j))$ 
 represents $\bbR \Hom^\bullet(\cO _{X/S},\cO _{X/S}(j))$. 

\begin{example} Let $X/S$ be a {\it Calabi-Yau} hyperplane configuration 
scheme, i.e. $d=n+2$. Then the poset $S$ is isomorphic to the poset of 
simplices in the standard triangulation of the boundary of an $(n+1)$-
simplex. Thus the complex $\Hom ^\bullet ({\check C}^\bullet, \cO _{X/S})$ 
is isomorphic to the (shifted by $n$ to the right) simplicial 
complex of the boundary of an $(n+1)$-simplex.
\end{example}

\subsection{Localization theorem}
Let $X_0$ be the reduced scheme, which is the union $X_0=L_1\cup ... 
\cup L_d\subset \bbP^{n+1}$. We have the obvious localization functor 
$$\cL :coh(X_0)\to coh(X/S).$$
Namely, for $\alpha \in S$ denote the closed embedding $t_{\alpha }:
\bbP _{\alpha }\to X_0$. Then given $F\in coh(X_0)$ we define 
$\cL (F)_{\alpha }:=t^*_{\alpha }F\in coh(\bbP _{\alpha })$. 
The functor $\cL$ is right exact, 
hence induces its left derived 
$$\bbL\cL :D^b(coh(X_0))\to D^-(coh(X/S)).$$

We will denote $\bbL\cL (\cO _{X_0}(i))$ simply by $\cO _{X/S}(i)$. 
The endofunctor 
$$\cO _{X/S}(i)\otimes _{\cO _{X/S}}(\cdot):coh(X/S)\to coh(X/S)$$
is an autoequivalence.

\begin{defn} A {\it perfect complex} in $D^b(coh(X_0))$ is a finite 
complex consisting of locally free sheaves. Denote by $\Perf (X_0)$ 
the full subcategory of $D^b(coh(X_0))$ consisting of perfect 
complexes. 
\end{defn}

Note that for a perfect complex $P^\bullet\in \Perf (X_0)$ we have 
$\cL (P^\bullet)=\bbL\cL(P^\bullet)\in D^b(coh(X/S))$. 

\begin{thm} The functor $\cL :\Perf (X_0)\to D^b(coh(X/S))$ is fully 
faithful.
\end{thm}

We need a few lemmas.

\begin{lemma}
For $i\geq 0$ the spaces $H^\bullet (X_0,\cO_{X_0}(i))$ and 
$H^\bullet (X/S,\cO _{X/S}(i))$ are isomorphic.
\end{lemma}

\begin{pf}
The scheme $X_0$ has a covering by closed subschemes 
$$X_0=\bigcup_{d(\alpha )=n}\bbP_{\alpha }.$$
Consider the corresponding \v{C}ech resolution 
$\cO _{X_0}(i)\to J^\bullet$ in $coh(X_0)$, where 
$$J^\bullet =\bigoplus_{d(\alpha)=n}\cO _{\bbP_{\alpha }}(i)\to 
\bigoplus_{d(\beta)=n-1}
\cO _{\bbP_{\beta }}(i)\to ...$$
(The fact that $\cO _{X_0}(i)\to J^\bullet$ is a quasiisomorphism 
can be easily checked by a local computation, since the hyperplanes $L_j$ 
are in general position). 
Note that each sheaf $\cO _{\bbP _{\gamma }}(i)$ is acyclic for 
$\Hom(\cO_{X_0}, 
\cdot )$ (since $i\geq 0$), 
hence the complex $\Hom ^\bullet(\cO _{X_0},J^\bullet)$ represents 
$\bbR \Hom ^\bullet(\cO_{X_0},\cO_{X_0}(i))$. 

On the other hand $\bbR \Hom ^\bullet(\cO _{X/S},\cO _{X/S}(i))$ is 
represented by the complex $\Hom ^\bullet ({\check C}^\bullet, 
\cO _{X/S}(i))$, which is naturally isomorphic to 
$\Hom ^\bullet(\cO _{X_0},J^\bullet)$. This proves the lemma.  
\end{pf}

\begin{cor}
We have $H^j(X/S,\cO _{X/S}(i))=0$ for $j>0$, $i>0$.
\end{cor}

\begin{pf} This follows from the last lemma and the vanishing 
of the corresponding cohomology groups on $X_0$. 
\end{pf}

\begin{lemma} For all $j\geq i$ the map 
$$\cL :\Hom (\cO _{X_0}(i),\cO _{X_0}(j))\to \Hom (\cO _{X/S}(i), 
\cO _{X/S}(j))$$
is an isomorphism.
\end{lemma}

\begin{pf} Since $\cL $ commutes with tensoring by $\cO _{X_0}(i)$ 
(resp. $\cO _{X/S}(i)$), it suffices to show that for all $j\geq 0$ the map 
$$\cL :H^0(X_0,\cO _{X_0}(j))\to H^0(X/S,\cO _{X/S}(j))$$
is an isomorphism. By Lemma 4.6 the two spaces have the same dimension. 
Also it is clear that the map is injective. Hence it is an isomorphism. 
\end{pf}

\begin{lemma} For any $F\in coh(X_0)$ and any $i$ the map 
$$\bbL\cL:\Ext^\bullet (F,\cO_{X_0}(i))\to \Ext^\bullet 
(\bbL\cL(F),\cO _{X/S}(i))$$
is an isomorphism. 
\end{lemma}

\begin{pf} We can find a left resolution $P^\bullet \to F$, where elements of 
$P^\bullet$ are direct sums of sheaves $\cO _{X_0}(-j_k)$ with $j_k>-i$. 
Then $\Hom ^\bullet (P^\bullet ,\cO _{X_0}(i))$ represents 
$\bbR \Hom ^\bullet (F,\cO _{X_0}(i))$ and by Corollary 4.7  
$\Hom ^\bullet (\cL (P^\bullet ),\cO _{X/S}(i))$ represents 
$\bbR \Hom ^\bullet (\bbL\cL (F),\cO _{X/S}(i))$. But the functor 
$\cL $ establishes an isomorphism of complexes  
$\Hom ^\bullet (P^\bullet , \cO_{X_0}(i))$ and 
$\Hom ^\bullet (\cL (P^\bullet), \cO_{X/S}(i))$ by Lemma 4.8
\end{pf}

\begin{pf} (of Theorem). It follows from the last lemma that the functor 
$\bbL\cL$ is fully faithful on the triangulated subcategory of 
$D^b(coh(X_0))$, 
generated by line bundles $\cO _{X_0}(i)$. To prove the theorem it 
suffices to show that any locally free sheaf $F$ on $X_0$ is a direct summand 
(in $D^b(coh(X_0))$) of a finite complex $P^\bullet $ which consists of 
direct sums of line bundles $\cO _{X_0}(i)$. For any $s\geq 0$ we can find an 
exact complex in $coh(X_0)$
$$0\to T\to P^{-s}\to ...\to P^{-1}\to P^0\to F\to 0,$$
where each $P^{-j}$ is a direct sum of line bundles $\cO _{X_0}(i)$. 
It represents an element in $\Ext ^{s+1}(F,T)$. Note that the functor 
$\Hom (F,(\cdot) )$ is isomorphic to $\Gamma (X_0,(\cdot)\otimes _{\cO _{X_0}}
F^*)$. Thus $\Ext^{s+1}(F,T)=0$ for $s>>0$ and so $P^\bullet \simeq F
\oplus T[s]$.
\end{pf}

\subsection{Exceptional collections and representations of quivers}
For a hyperplane configuration scheme $X/S$ we will prove that the derived 
category $D^b(coh(X/S))$ is equivalent to the derived category 
of finite dimensional modules over a certain finite dimensional algebra $A$. 
Let us first recall some basic definitions (see, for example, [Bo]). 

Let $\cA$ be a triangulated category.

\begin{defn} An object $E$ in $\cA$ is called {\it exceptional} if 
$\Ext ^i(E,E)=0$ for $i\neq 0$, $\Hom (E,E)=\bbC$. 
\end{defn}

\begin{defn} An ordered collection $(E_0,...,E_m)$ of exceptional objects  
is called an {\it exceptional collection} in $\cA$ if $\Ext^i(E_j,E_k)=0$ 
for all $i$ when $j>k$. 
\end{defn} 

\begin{defn} An exceptional collection $(E_0,...,E_m)$ in $\cA$ satisfying 
the condition $\Ext ^i(E_j,E_k)=0$ for all $j$ and $k$ if $i\neq 0$ is called 
a {\it strong exceptional collection}. 
\end{defn}

Given a strong exceptional collection $(E_0,...,E_m)$ in $\cA$ put 
$E=\oplus E_i$ and $A=\Hom(E,E)$. Then $A$ is a finite dimensional 
algebra. It is actually the {\it path algebra} of a certain quiver 
with relations which contains $n+1$ vertices (see [Bo] for the 
discussion of this topic). The following theorem is proved in [Bo].

\begin{thm}
Let $Y$ be a smooth projective variety. Assume that the derived category 
$\cA =D^b(coh(Y))$ is generated (as a triangulated category) by a strong 
exceptional collection $(E_0,...,E_m)$. Then $\cA $ is equivalent to the 
derived category $D^b(A-mod)$ of right finite-dimensional $A$-modules. 
\end{thm}

\begin{remark}
Strong exceptional collections as in the above theorem are known 
for projective spaces,  quadrics, Grassmannians and flag varieties.
\end{remark}

\begin{remark}
Let $p_k:E\to E_k$ be the projection on the $k$-th summand. Consider the 
projective right $A$-module $P_k:=p_kA$. We have the natural isomorphism 
$$\Hom _A(P_k,P_l)=\Hom (E_k,E_l).$$
\end{remark}

Here we prove a counterpart of the above theorem for hyperplane 
configuration schemes. Namely, let $X/S$ be a hyperplane configuration 
scheme of type $(n,d)$. For $\alpha \in S$ put  
$$E_{\alpha ,i}:=\cO _{U_{\alpha }}(i)\in coh(X/S).$$
Consider the collection 
$$\cE =\{ E_{\alpha ,i}\vert \alpha \in S, 0\leq i \leq dim\bbP_{\alpha }\}.$$
It has a natural lexicographical partial order, i.e. $E_{\alpha ,i}
>E_{\beta ,j}$ iff $\alpha >\beta$, or $\alpha =\beta $ and $i>j$. 
Again we put 
$$E=\oplus_{E_{\alpha ,i}\in \cE}E_{\alpha , i},\quad A=\Hom(E,E).$$

\begin{thm}
Let $X/S$ be a hyperplane configuration scheme. Put $\cA =D^b(coh(X/S))$. 
Consider the collection $\cE $ as above with any linear ordering which is 
compatible with its natural partial order. Then 

i) $\cE $ is a strong exceptional collection in $\cA$;

ii) $\cE$ generates the triangulated category $\cA$;

iii) the category $\cA$ is equivalent to the derived category 
$D^b(A-mod)$ of finite dimensional right $A$-modules. 
\end{thm}

\begin{pf} i). Let $E_{\alpha ,i},E_{\beta ,j}\in \cE$. It follows from Lemma 
3.21 that 
$$\Ext ^\bullet (E_{\alpha ,i},E_{\beta ,j})=
H ^\bullet (\bbP_{\alpha}, (E_{\beta ,j-i})_{\alpha }).$$
Thus $\Ext^\bullet(E_{\alpha ,i},E_{\beta ,j})=0$ unless $\beta \geq \alpha$. 
Assume that $\beta \geq \alpha $. Then by the same formula 
$$ \Ext ^\bullet (E_{\alpha ,i},E_{\beta ,j})=
H ^\bullet (\bbP_{\alpha}, \cO _{\bbP_{\alpha }}(j-i)).$$
If $j<i$, then again $\Ext ^\bullet=0$ 
(because $i-j\leq dim(\bbP_{\alpha} )$). 
If $j\geq  i$ then only $\Ext ^0$ is nonzero. Thus $\cE$ is a strong 
exceptional collection.

ii). We will use the following theorem of A.Beilinson: for a projective space 
$\bbP^k$ the line bundles $\cO _{\bbP^k},\cO _{\bbP^k}(1),...,
\cO _{\bbP^k}(k)$ generate the category $D^b(coh(\bbP ^k))$. 

For $\alpha \in S$ denote by $i_{\alpha }:\{\alpha \}\hookrightarrow 
U_{\alpha }$ and $j_{\alpha }:U_{\alpha }\hookrightarrow S$ the closed 
and open embeddings respectively. Recall that the functor 
$$j_{\alpha !}i_{\alpha *}:D^b(coh(\bbP_{\alpha }))\to D^b(coh(X/S))$$
is fully faithful (Corollary 3.7). The images of these functors generate the 
triangulated category $D^b(coh(X/S))$. So it suffices to prove that the 
collection $\cE$ generates each subcategory $j_{\alpha !}i_{\alpha *}
D^b(coh(\bbP_{\alpha }))$. We will prove this assertion by induction on 
$\alpha $ (using the partial order of $S$).

If $\alpha $ is a smallest element the assertion follows from the 
above mentioned theorem of Beilinson about $D^b(coh(\bbP^k))$. 
Supposed we proved the 
assertion for each $\beta \in S$ such that $\beta <\alpha $. 
Put $V=U_{\alpha }-\{\alpha \}$. For $E_{\alpha ,i}\in \cE$ consider the 
exact sequence 
$$0\to (E_{\alpha ,i})_V\to E_{\alpha ,i}\to j_{\alpha !}i_{\alpha *}
\cO _{\bbP_{\alpha }}(i)\to 0.$$
By our assumption $(E_{\alpha ,i})_V$ is contained in the subcategory 
of $D^b(coh(X/S))$ generated by $\cE$. So the same is true for 
$j_{\alpha !}i_{\alpha *}
\cO _{\bbP_{\alpha }}(i).$ Now the assertion for $\alpha $ follows from the 
same fact about $\bbP^k$ mentioned in the beginning.

iii) The proof is essentially the same as that of the corresponding 
Theorem 6.2 in [Bo]. We present it for completeness.

Let $A-mod$ denote the abelian category of right $A$-modules, 
$(A-mod)_f\subset A-mod$ - the full subcategory of finitely 
generated $A$-modules. Let $D^b(A-mod)$, $D^b((A-mod)_f)$ be the 
corresponding bounded derived categories. Denote by 
$D^b_f(A-mod)\subset D^b(A-mod)$ the full subcategory consisting 
of complexes with cohomologies in $(A-mod)_f$. The natural 
functor $D^b((A-mod)_f)\to D^b_f(A-mod)$ is an equivalence. 
We are going to construct an equivalence 
$$\theta :D^b(coh(X/S))\stackrel{\sim}{\rightarrow}D_f^b(A-mod).$$

For each $F\in D^b(coh(X/S))$ choose an injective resolution 
$F\to I^\bullet(F)$. Put $\theta (F):= \Hom^\bullet(E,I^\bullet(F))$. 
The algebra $A$ acts on this complex from the right 
preserving each component and 
commuting with the differential. Thus $\theta (F)$ is a complex 
of right $A$-modules. The cohomologies $H^i(\theta (F))$ are finite 
dimensional vector spaces and are zero for $|i|>>0$. Thus $\theta (F)$ 
is an object in $D^b_f(A-mod)$ and we obtain a functor 
$$\theta :D^b(coh(X/S))\to D^b_f(A-mod).$$ 
Let us prove that  $\theta$ is an equivalence.

Denote by $p_{\alpha ,i}:E\to E_{\alpha ,i}$ the projection. Let 
$P_{\alpha ,i}:=p_{\alpha ,i}A$ be the corresponding projective 
right $A$-module.

Note that $H^j(\theta (E_{\alpha ,i}))=0$ for $j\neq 0$ and 
$H^0(\theta (E_{\alpha ,i}))=P_{\alpha ,i}$. The collections 
$\cE =\{E_{\alpha , i}\}$ and $P_{\alpha ,i}$ generate the categories 
$D^b(coh(X/S))$ and $D^b_f(A-mod)$ respectively. Thus it suffices to prove 
that the map 
$$\theta :\Ext ^\bullet (E_{\alpha ,i},E_{\beta ,j})\to 
\Ext ^\bullet (P_{\alpha ,i},P_{\beta ,j})$$
is an isomorphism. But the higher Ext's on both sides vanish and 
$$\theta :\Hom (E_{\alpha ,i},E_{\beta ,j})\to 
\Hom _A(P_{\alpha ,i},P_{\beta ,j})$$
is an isomorphism by Remark 4.15 above.
\end{pf}

\subsection{Spherical objects and twist functors}
We recall some facts about spherical objects and the corresponding 
twist functors following [ST]. In the next subsection we will 
apply these results to Calabi-Yau hyperplane configuration schemes.

Let $\cC$ be an abelian category and $\cC^\prime\subset \cC$ be a full
 abelian subcategory with the following properties

(C1) $\cC^\prime $ is a Serre subcategory of $\cC$;

(C2) $\cC$ contains infinite direct sums and products;

(C3) $\cC$ has enough injectives, and any direct sum of injectives is again 
an injective;

(C4) for any epimorphism $f:A\twoheadrightarrow A^\prime $ with $A\in \cC$ 
and $A^\prime \in \cC^\prime$, there is a $B^\prime \in \cC^\prime$ and 
$g:B^\prime \to A$ such that $fg$ is an epimorphism.

\begin{lemma}
Let $X$ be a Noetherian scheme over a field $k$ and $\cC=qcoh(X)$, 
$\cC^\prime =coh(X)$. Then properties (C1)-(C4) are satisfied. 
\end{lemma}

\begin{pf} See [ST].
\end{pf}

\begin{defn}
Let $K^+(\cC)$ be the homotopy category of bounded below complexes of 
objects in $\cC$. Denote by $\cK\subset K^+(\cC)$ the full subcategory 
whose objects are complexes $C^\bullet$ of injectives in $\cC$ which 
satisfy $H^i(C^\bullet)\in \cC^\prime$ for all $i$ and $H^i(C^\bullet)=0$ 
for all $i>>0$. 
\end{defn}

The following lemma is well known (also proved in [ST]).

\begin{lemma}
The obvious functors 
$$D^b(\cC^\prime)\rightarrow D^b_{\cC^\prime}(\cC)\leftarrow \cK$$ 
are equivalences. Thus in particular $D^b(\cC^\prime)\simeq \cK$.
\end{lemma}

Let us define some operations in the category $\cK$.

\begin{defn} Let $E\in \cK$ be an object satisfying the following finiteness 
conditions

(K1) $E$ is a bounded complex,

(K2) for any $F\in \cK$ 
$$\bigoplus_idim\Ext^i(E,F)<\infty,\quad \bigoplus_idim\Ext^i(F,E)<\infty.$$
Then one defines the {\it twist functor} $T_E:\cK \to \cK$ by 
$$T_E(F):=Cone(\Hom^\bullet(E,F)\otimes E\stackrel{ev}{\rightarrow}F),$$
where $ev$ is the obvious evaluation map. 
\end{defn}
\medskip
\noindent{\bf Remarks.}

\noindent{1.} $T_E$ is an exact functor , i.e. it takes exact triangles to 
exact triangles.

\noindent{2.} If $E, E_1 \in \cK$ satisfying (K1),(K2) are isomorphic, 
then the corresponding functors $T_E,$ $T_{E_1}$ are isomorphic.
\medskip
\begin{defn} An object $E\in \cK$ is called $n$-{\it spherical} for some 
$n>0$ if it satisfies (K1),(K2) above and in addition 

(K3) $\Ext ^i(E,E)$ is equal to $\bbC$ for $i=0,n$ and zero otherwise;

(K4) the composition 
$$\Ext ^j(E,F)\times \Ext ^{n-j}(F,E)\to \Ext ^n(E,E)\simeq \bbC$$
is a nondegenerate pairing for all $F\in \cK$, $j\in \bbZ$.
\end{defn}

\begin{prop} If $E$ is $n$-spherical for some $n>0$, then the functor 
$T_E:\cK \to \cK$ is an equivalence.
\end{prop}

\begin{pf} See [ST].
\end{pf}
   
Using the equivalence $D^b(\cC^\prime)\simeq \cK$ we obtain the following 
corollary.

\begin{cor} In the setup of the last proposition the $n$-spherical 
object $E$ induces an exact auto-equivalence of $D^b(\cC^\prime)$ 
which we also denote $T_E$. 
\end{cor}

\begin{example} The standard example of an $n$-spherical object 
in $D^b(coh(X))$ is $\cO _X$ for a smooth $n$-dimensional Calabi-Yau 
variety $X$. Indeed, by definition 
$\Ext ^i(\cO_X,\cO_X)$ is $\bbC$ for $i=0,n$ and zero otherwise.
Also by definition $\cO _X\simeq \omega _X$ and 
hence by Serre duality the natural pairing 
$$\Ext^i(\cO _X,F)\times \Ext ^{n-i}(F,\cO _X)\to \Ext ^n(\cO_X,\cO_X)
\simeq \bbC$$
is perfect for all $F\in D^b(coh(X))$ and all $i$. Thus $\cO_X$ is 
indeed $n$-spherical. Examples of other spherical objects in $D^b(coh(X))$ 
are discussed in [ST]. 
\end{example}

\subsection{Calabi-Yau hyperplane configuration schemes}
Let $X/S$ be a Calabi-Yau hyperplane configuration scheme, i.e. it is 
of type $(n,n+2)$ for some $n>0$. Let $X_0$, as usual, be the corresponding 
singular subscheme of $\bbP ^{n+1}$ which is the union of the $n+2$ 
hyperplanes. The object $\cO _{X_0}\in D^b(coh(X_0))$ is the dualizing 
complex on $X_0$. That is the pairing 
$$\Ext ^i(\cO_{X_0},F)\times \Ext ^{n-i}(F,\cO_{X_0})\to \Ext ^n
(\cO_{X_0},\cO_{X_0})\simeq \bbC$$
is perfect for all $F\in D^b(coh(X_0))$ and all $i$. 
Thus for all $F\in coh(X_0)$ we have $\Ext ^{n+1}(F,\cO _{X_0})=0$. 
Therefore the same is true for all $F\in qcoh(X_0)$. 
It follows that $\cO_{X_0}$ has finite 
injective dimension in $qcoh(X_0)$. 
So there exists a finite injective resolution 
$\cO_{X_0}\to E$ in $qcoh(X_0)$ such that $E\in \cK$ satisfies the 
properties (K1)-(K4) above. Hence the corresponding twist 
functor $T_{\cO_{X_0}}=T_E:D^b(coh(X_0))\to D^b(coh(X_0))$ is an equivalence. 
It follows that $\cO_{X_0}(j)$ is also $n$-spherical for all $j\in \bbZ$. 
Hence for each $j$ we obtain the corresponding auto-equivalence of 
$D^b(coh(X_0))$. We are going to show that these auto-equivalences 
survive when we pass from $X_0$ to the configuration scheme $X/S$. 

First we need a lemma.

\begin{lemma} The categories  $\cC=qcoh(X/S)$ and $\cC^\prime =coh(X/S)$ 
satisfy the assumptions (C1)-(C4) of the last subsection.
\end{lemma}

\begin{pf} (C1) is obvious and (C4) was proved in Lemma 3.4 above. 
Let us prove (C2). 

Let $F^i\in qcoh(X/S)$ for $i\in I$. Let us define the object 
$F=\oplus F^i\in qcoh(X/S)$. Namely, for $\alpha \in S$ put $F_{\alpha }=
\oplus F^i_{\alpha }$. Then for $f_{\beta \alpha}:X(\beta )\to X(\alpha )$ 
we have 
$$f^*_{\beta \alpha}(\bigoplus F^i_{\alpha})=
\bigoplus f^*_{\beta \alpha }F^i_{\alpha}.$$
Hence we define the morphism $\sigma _{\alpha \beta}(F):f^*_{\beta \alpha}
F_{\alpha }\to F_{\beta }$ to be the direct sum of morphisms $\sigma 
_{\alpha \beta}(F^i)$. 

The definition of $G=\prod F^i\in qcoh(X/S)$ is similar. Namely, 
for $\alpha \in S$ put $G_{\alpha }=\prod F^i_{\alpha }$. Then for 
$f_{\beta \alpha}:X(\beta )\to X(\alpha )$ we have a canonical morphism
$$m_{\alpha \beta}:f^*_{\beta \alpha}(\prod F^i_{\alpha })\to 
\prod f^*_{\beta \alpha}F^i_{\alpha}.$$
Hence we define the morphism $\sigma _{\alpha \beta}(G):f^*_{\beta \alpha}
G_{\alpha }\to G_{\beta}$ as the composition of $m_{\beta \alpha}$ with 
the product of morphisms $\sigma _{\alpha \beta}(F^i)$. 

One checks directly that $F$ and $G$ satisfy the universal properties 
of the direct sum and the direct product respectively. 

Let us prove (C3). 
Also by Lemma 3.1 the category $qcoh(X/S)$ has enough injectives. 
It remains to prove that a direct sum of injectives in $qcoh(X/S)$ is 
again an injective. Note that the category $qcoh(X/S)$ is locally 
noetherian, i.e. it has a family of noetherian generators (any object 
in $coh(X/S)$ is noetherian). Hence a direct sum of injectives in $qcoh(X/S)$ 
is an injective by Corollary 6.50 in [BD].
\end{pf}

\begin{remark}
The assertion of the lemma is true for arbitrary configuration schemes (the 
proof is the same).
\end{remark}

\begin{thm} Let $X/S$ be a Calabi-Yau hyperplane configuration 
scheme of type $(n,n+2)$. Then for any $t$ the line bundle $\cO _{X/S}(t)$ is  
 an $n$-spherical object in $D^b(coh(X/S))$.
\end{thm}

\begin{pf} We may assume that $t=0$. 

The finiteness conditions (K1),(K2) are satisfied  for any 
$E\in D^b(coh(X/S))$ by Corollary 3.2 and Lemma 3.12. Also by Theorem 
4.5 $\Ext^i(\cO _{X/S},\cO _{X/S})=\Ext ^i(\cO _{X_0},\cO _{X_0})$. 
So we only need to prove that for all $F\in D^b(coh(X/S))$ 
the natural pairing 
$$\Ext^i(\cO _{X/S}, F)\times \Ext ^{n-i}(F,\cO_{X/S})\to 
\Ext ^n(\cO _{X/S},\cO _{X/S})\simeq \bbC$$
is perfect.

\medskip
\noindent{\it Step 1. } We can find a left resolution $P^\bullet 
\to F$, where $P^\bullet$ consists of finite direct sums of sheaves 
$\cO _{U_{\alpha }}(j)$ for various $\alpha \in S$, $j<0$. Thus we 
may assume that $F =\cO _{U_{\alpha }}(j)$ for some $\alpha \in S$ 
and $j<0$. 

\medskip
\noindent{\it Step 2.} Note that 
$$\Ext ^k(\cO _{U_{\alpha }}(j),\cO _{X/S})=\Ext ^k(\cO _{U _{\alpha }}(j), 
\cO _{U_{\alpha }})=\Ext ^k(\cO _{\bbP _{\alpha }}(j), 
\cO _{\bbP _{\alpha }})=0,$$
unless $k=0$. 

\begin{lemma} $\Ext ^i(\cO _{X/S},\cO _{U_{\alpha }}(j))=0$ if $i\neq n$.
\end{lemma}

\noindent{\it Proof of lemma.} Consider the \v{C}ech 
resolution $\check C^\bullet \to \cO _{X/S}$ as defined above.

By Lemma 3.21 we have 
$$\Ext ^p(\cO _{U_{\gamma }},\cO _{U_{\alpha }}(j))= 
\begin{cases}
             0, &\ \text{if $\gamma \nleq \alpha$ }, \\
 \Ext ^p(\cO _{\bbP _{\gamma }},\cO _{\bbP _{\gamma }}(j)),  
&\ \text{otherwise}.
\end{cases}
$$
%$$(F_U)_{\alpha }=\begin{cases}
%F_{\alpha}, &\ \text{if $\alpha\in U$}, \\
%0,          &\ \text{if $\alpha\in Z$}.
%\end{cases} 
%$$
Since $j<0$, $\Ext ^p(\cO _{\bbP _{\gamma }}, \cO _{\bbP _{\gamma }}(j))=0$, 
unless $p=d(\alpha )$. So the spectral sequence argument applied 
to the naturally filtered complex
$\bbR\Hom ^\bullet(\check C^\bullet,\cO _{U_{\alpha }}(j))$ yields 
$\Ext ^i(\cO _{X/S},\cO _{U_{\alpha }}(j))=0$ for $i\neq n$. This 
proves the lemma. 

\medskip
Therefore, it suffices to show that the pairing 
$$\leqno (1) \Ext ^n(\cO _{X/S},\cO _{U_{\alpha }}(j))\times 
\Hom(\cO _{U_{\alpha }}(j),\cO _{X/S})\to 
\Ext ^n(\cO _{X/S},\cO _{X/S})$$
is perfect.

\medskip
\noindent{\it Step 3.} 
Let $\check C_{\alpha }^\bullet\subset \check C^\bullet$ be the 
subcomplex consisting of direct sums of objects $\cO _{U_{\beta }}$ 
with $\beta \leq \alpha$. 

\begin{lemma} The inclusion $\check C^\bullet_{\alpha}\hookrightarrow 
\check C^\bullet$ induces isomorphisms

 $$\Ext ^n(\check C^\bullet,\cO _{U_{\alpha }}(j))\stackrel{\sim}{\to}
\Ext^n(\check C^\bullet_{\alpha },\cO _{U_{\alpha }}(j)),$$
and 
 $$\Ext ^n(\check C^\bullet,\cO _{X/S})\stackrel{\sim}{\to}
\Ext^n(\check C^\bullet_{\alpha},\cO _{X/S}).$$
\end{lemma}

\noindent{\it Proof of lemma.} The first one follows from Lemma 3.21. 
Let us prove the second one.
 The complex $\Hom ^\bullet (\check C^\bullet, \cO _{X/S})$ is 
isomorphic to the (shifted by $n$ to the right) 
chain complex of the boundary of an $(n+1)$-
simplex with the standard triangulation. Its quotient $\Hom ^\bullet
(\check C^\bullet_{\alpha},\cO _{X/S})$  
 corresponds to the (shifted) chain complex of a subsimplex. Hence 
the projection $\Hom ^\bullet (\check C^\bullet, \cO _{X/S})\to  
\Hom ^\bullet
(\check C^\bullet_{\alpha},\cO _{X/S})$  
induces an isomorphism on the top cohomology.
This proves the lemma. 

\medskip
It follows that the pairing (1) is isomorphic to 
$$\Ext ^n(\check C^\bullet _{\alpha },\cO _{U_{\alpha }}(j))\times 
\Hom (\cO _{U_{\alpha }}(j), \cO _{X/S})
\to \Ext ^n(\check C^\bullet _{\alpha },\cO _{X/S}).$$
Now again by Lemma 3.21 we can replace $\cO _{X/S}$ by 
$\cO _{U_{\alpha }}$ to get 
the isomorphic pairing
$$\Ext ^n(\check C^\bullet _{\alpha },\cO _{U_{\alpha }}(j))\times 
\Hom (\cO _{U_{\alpha }}(j), \cO _{U_{\alpha}})
\to \Ext ^n(\check C^\bullet _{\alpha },\cO _{U_{\alpha }}).$$
 
Let $j_{\alpha }:\{\alpha \}\hookrightarrow U_{\alpha }$ be the closed 
embedding. Note that $\cO _{U_{\alpha }}=j_{\alpha +}\cO _{\bbP _{\alpha}}$, 
and $\cO _{U_{\alpha }}(j)=j_{\alpha +}\cO _{\bbP _{\alpha}}(j)$. Thus the 
last pairing is 

$$\leqno (2) \Ext ^n(\check C^\bullet _{\alpha },j_{\alpha +}
\cO _{\bbP _{\alpha}}(j))
\times 
\Hom (j_{\alpha +}\cO _{\bbP _{\alpha}}(j),j_{\alpha +}\cO _{\bbP _{\alpha}})
\to \Ext ^n(\check C^\bullet _{\alpha },j_{\alpha +}\cO _{\bbP _{\alpha}}).$$

\noindent{\it Step 4.}
\begin{lemma} Let $F\in D^b(coh(\bbP _{\alpha }))$ be a finite complex 
of locally free sheaves. Then for $i>>0$ we have a natural quasi-isomorphism 
of complexes
$$\Hom ^\bullet
(\cO _{\bbP _{\alpha }}(-i),F(-d(\alpha)-1))
\stackrel{\sim}{\to } \Hom ^\bullet
(\check C^\bullet _{\alpha }(-i) , 
 j_{\alpha +}F)[n-d(\alpha)],$$
where the first $\Hom$ is computed in the category 
$coh(\bbP _{\alpha})$ and the second in $coh(X/S)$ (or in 
$coh(X/U_{\alpha })$).
\end{lemma}

\noindent{\it Proof of lemma.} 
Let $\beta _1,...,\beta _{d(\alpha )+1}\in S$ be all elements with the 
property $\beta _i<\alpha$, $d(\beta _i)=d(\alpha )-1$. Then for each $i$,  
$\bbP _{\beta _{i}}$ is a hyperplane in $\bbP _{\alpha }$ (more precisely, 
$f_{\beta _i\alpha}(\bbP _{\beta _i})$ is such). Put $Y=\cup \bbP _{\beta _i}
\subset \bbP _{\alpha }$.  
We have the exact sequence of complexes  
on $\bbP _{\alpha }$
$$0\to F(-d(\alpha )-1)\to F\to F_Y\to 0.$$
Let $0\to F_Y\to J^\bullet (F_Y)\to 0$ be the \v{C}ech resolution (as in 
the proof of Lemma 4.6) of 
$F_Y$ corresponding to the closed covering of $Y$ by 
$\bbP _{\beta _i}$'s, i.e. 
$$J^\bullet (F_Y)=\ 0\to \bigoplus _iF_{\bbP _{\beta _i}}\to 
\bigoplus _{i<k}F_{\bbP _{\beta _i}\cap \bbP _{\beta _j}}\to ...$$
Denote by $J^\bullet _F$ the complex $0\to F \to J^\bullet(F_Y)\to 0$. 
Then  $J^\bullet _F$ is quasi-isomorphic to the complex $F(-d(\alpha )-1)$. 
Thus for $i>>0$ the natural map of complexes 
$$\Hom^\bullet(\cO _{\bbP _{\alpha }}(-i),F(-d(\alpha)-1))\to 
\Hom^\bullet(\cO _{\bbP _{\alpha }}(-i),J^\bullet _F)$$
is a quasi-isomorphism. It remains to notice the natural isomorphism 
of complexes
$$\Hom ^\bullet_{coh(\bbP _{\alpha })}(\cO _{\bbP _{\alpha }}(-i), 
J^\bullet _F)
= \Hom ^\bullet _{coh(X/S)}(\check C^\bullet _{\alpha }(-i), 
j_{\alpha +}F)[n-d(\alpha )].$$
This proves the lemma.

\medskip
\noindent{\it Step 5.} Let $F\in D^b(coh(\bbP _{\alpha}))$ be as in the 
above lemma. We may take a left resolution 
$P^\bullet \to \cO _{\bbP _{\alpha}}$ consisting of finite sums 
$\oplus \cO _{\bbP _{\alpha}}(-i)$ for $i>>0$. This induces a 
quasi-isomorphism of \v Cech resolutions 
$\check C^\bullet _{\alpha}(P^\bullet)\stackrel{\sim}{\rightarrow}
\check C^\bullet_{\alpha}$. We have 
$$\Hom^\bullet(P^\bullet,F(-d(\alpha)-1))=\bbR \Hom^\bullet
(\cO _{\bbP _{\alpha}},F(-d(\alpha)-1)),$$
$$\Hom^\bullet(\check C^\bullet _{\alpha}(P^\bullet),j_{\alpha +}(F))
=\bbR \Hom ^\bullet(\check C^\bullet _{\alpha},j_{\alpha +}(F)).$$
Thus the above lemma implies a natural quasi-isomorphism of 
complexes
$$\bbR \Hom^\bullet(\cO _{\bbP _{\alpha }},F(-d(\alpha)-1))
\stackrel{\sim}{\rightarrow}\bbR \Hom ^\bullet (\check C^\bullet
_{\alpha},j_{\alpha +}F)[n-d(\alpha)].$$

This shows that the pairing (2) above is isomorphic to the pairing 
$$\Ext ^{d(\alpha )}(\cO _{\bbP _{\alpha }},\cO _{\bbP _{\alpha }}
(j-d(\alpha )-1))\times \Hom (\cO _{\bbP _{\alpha }}(j-d(\alpha )-1), 
\cO _{\bbP _{\alpha }}(-d(\alpha) -1))$$
$$\to 
\Ext ^{d(\alpha )}(\cO _{\bbP _{\alpha}},\cO _{\bbP _{\alpha }}
(-d(\alpha )-1)),$$
where everything is computed in $coh(\bbP _{\alpha })$. The   
last pairing is perfect by the Serre duality on $\bbP _{\alpha }$. 
This proves the theorem.
\end{pf}

\begin{cor} Let $X/S$ be a Calabi-Yau hyperplane configuration 
scheme of type $(n,n+2)$. Then for any $t\in \bbZ$ the twist functor 
$T_{\cO_{X/S}}(t)$ is an autoequivalence of $D^b(coh(X/S)$.
\end{cor}

\begin{pf} This follows immediately from the last theorem, Lemma 4.25
 and Corollary 4.23.
\end{pf}

\end{document}